\documentclass[10pt,twoside]{amsart}
\usepackage{amsmath,amsfonts,amssymb,mathrsfs,amsthm,charter,graphicx}



\newcommand{\mat}{\mathbb{M}(n)}
\newcommand{\De}{{\rm D}}
\newcommand{\tr}{\mathop{{\rm tr}}}
\newcommand{\padj}{\mathop{{\rm padj}}}
\newcommand{\adj}{\mathop{{\rm adj}}}
\newcommand{\per}{\mathop{{\rm per}}}
\newcommand{\I}{\mathcal{I}}
\newcommand{\J}{\mathcal{J}}
\newcommand{\K}{\mathcal{K}}
\newcommand{\Hil}{\mathcal{H}}
\newcommand{\lo}{\mathcal{B}}

\newcommand{\N}{\mathbb N}
\newtheorem{theorem}{Theorem}

\newtheorem{remark}[theorem]{Remark}

\title{
Derivatives of tensor powers and their norms}

\author{Rajendra Bhatia}
\address{ Theoretical Statistics and Mathematics Unit, Indian Statistical Institute, 7, S.J.S. Sansanwal Marg, New Delhi-110016, India.}
\author{Priyanka Grover}
\address{ Theoretical Statistics and Mathematics Unit, Indian Statistical Institute, 7, S.J.S. Sansanwal Marg, New Delhi-110016, India.}
\author{Tanvi Jain}
\address{ Theoretical Statistics and Mathematics Unit, Indian Statistical Institute, 7, S.J.S. Sansanwal Marg, New Delhi-110016, India.}
\email{rbh@isid.ac.in, tanvi@isid.ac.in, pgrover8r@isid.ac.in}
\begin{document}
\bibliographystyle{plain}



\pagestyle{myheadings}
\markboth{R.\ Bhatia, P.\ Grover, and T.\ Jain}{Derivatives of tensor powers and their norms}
\maketitle

 \begin{abstract}
The norm of the $m$th derivative of the map that takes an operator to its $k$th antisymmetric tensor power is evaluated. The case $m=1$ has been studied earlier by  Bhatia and  Friedland [R.~Bhatia and S.~Friedland.
\newblock Variation of Grassman powers and spectra. 
\newblock {\em Linear Algebra and its Applications}, 40:1--18, 1981]. For this purpose a multilinear version of a theorem of Russo and Dye is proved: it is shown that a positive $m$-linear map between $C^{\ast}$-algebras attains its norm at the $m$-tuple $(I, \, I, \ldots, I).$ Expressions for derivatives of the maps that take an operator to its $k$th tensor power and $k$th symmetric tensor power are also obtained. The norms of these derivatives are computed. Derivatives of the map taking a matrix to its permanent are also evaluated.
 \end{abstract}

\numberwithin{theorem}{section}
\numberwithin{equation}{section}

\section{Introduction}
 Let $\mathcal{B}(\mathcal{H})$ be the space of linear operators on an $n$ dimensional Hilbert space $\mathcal{H}.$ Let $s_1 (A) \ge s_2 (A) \ge \cdots \ge s_n (A) \ge 0$ be the decreasingly ordered singular values of an operator $A.$ Let $\wedge^k \Hil$ be the $k$th antisymmetric tensor power of $\Hil,\ 1\leq k\leq n;$ and let $\wedge^k: \mathcal B(\Hil)\rightarrow \mathcal B(\wedge^k \Hil)$ be the map that takes an element $A$ of $\Hil$ to its $k$th antisymmetric tensor power $\wedge^k(A).$ Let $\De \wedge^k(A)$ be the derivative of this map. This is a linear map from $\mathcal B(\Hil)$ into $\mathcal B(\wedge^k \Hil)$ and its norm is defined as 
\begin{equation}
\|\De \wedge^k (A)\|=\underset{\|X\|=1}{\sup}\ \|\De \wedge^k(A)(X)\|, \label{1}
\end{equation}
where $\|X\|$ is the operator norm of $X$ as a linear operator on $\Hil;$ i.e.,
\begin{equation}
\|X\|=\underset{\|u\|=1}{\sup}\ \|X u\|.\label{2}
\end{equation}
 An alternative expression for this is $\|X\|=s_1(X).$ Finding $\|X\|$ involves solving a maximisation problem, which is not easy.

Motivated by problems in perturbation theory of eigenvalues, R. Bhatia and S. Friedland \cite{rbhfriedland} studied the problem of finding the norm \eqref{1} and obtained a striking formula

\begin{equation}
\|\De \wedge^k(A)\|=p_{k-1}(s_1(A),\ldots,s_k(A)), \label{3}
\end{equation}
where $p_{k-1}(x_1,\ldots,x_k)$ is the $(k-1)$th elementary symmetric polynomial in $k$ variables $x_1,\ldots,x_k.$
Analogues of this formula for other kinds of tensor products have been established; see \cite{rbh}, \cite{rbhsilva}.

This paper is concerned with higher order derivatives of $\wedge^k$ and of other multilinear operators and functions, and is related to some other recent work of two of the authors. The famous Jacobi formula gives the derivative of the determinant function on $n\times n$ matrices as 
\begin{equation}
\De \det (A)(X)=\tr (\adj (A)X),\label{Jacobi}
\end{equation}
where the symbol $\adj{(A)}$ stands for the transpose of the matrix whose $(i,j)$-entry is $(-1)^{i+j} \det {A(i|j)}$, called the \textit{adjugate} (or the \textit{classical adjoint}) of $A$. Analogous formulas for higher order derivatives of $\det$ were obtained in \cite{rbhjain}. Then the more general problem of evaluating higher order derivatives of the map $\wedge^k$ was studied in \cite{tjain}. For $1\leq m \leq k,$ $\De^m \wedge^k (A)$ is a multilinear map
\begin{equation}
\De^m \wedge^k (A) : \underbrace{\mathcal B(\Hil) \times \cdots \times \mathcal B(\Hil)}_{m\text{-fold}} \rightarrow \mathcal B(\wedge^k \Hil). \label{5'}
\end{equation}
Its norm is defined as
\begin{equation}
\|\De^m \wedge^k (A)\|=\underset{\|X^1\|=\cdots=\|X^m\|=1} {\sup}\ \|\De^m \wedge^k (A) (X^1,\ldots,X^m)\|.\label{6'}
\end{equation}

In \cite{tjain} Jain obtained a formula for \eqref{5'} and used it to evaluate \eqref{6'}. This last result can be stated as:

\begin{theorem}\label{theorem1}
With notations as above, we have
\begin{equation}
\|\De^m \wedge^k (A)\|=m!\ p_{k-m} (s_1(A),\ldots,s_k(A)), \label{7}
\end{equation}
where $p_{k-m}$ is the $(k-m)$th elementary symmetric polynomial.
\end{theorem}

Recall that for $1\leq r\leq k,$ the $r$th elementary symmetric polynomial is defined as 
\begin{equation*}
p_r(x_1,\ldots,x_k)=\sum_{1\leq i_1<\cdots<i_r\leq k} x_{i_1} x_{i_2}\cdots x_{i_r} .
\end{equation*} 

The first step in the Bhatia-Friedland proof of \eqref{3} is the observation that
\begin{equation}
\|\De \wedge^k (A)\|=\|\De \wedge^k (|A|)\|, \label{8}
\end{equation}
where $|A|$ is the absolute value of $A,$ defined as $|A|=(A^*A)^{1/2}.$ This was exploited by V. S. Sunder \cite{sunder}, who obtained another proof of \eqref{3} by invoking a well-known theorem from the theory of positive linear maps. A linear map $\Phi$ from $\mathcal B(\Hil)$ into $\mathcal B(\K)$ is said to be \emph{positive} if $\Phi(A)$ is a positive semidefinite operator whenever $A$ is positive semidefinite. A famous theorem of Russo and Dye \cite[p.42]{rBh} says that if $\Phi$ is a positive linear map, then $\|\Phi\|=\|\Phi(I)\|.$ The proof that we give for \eqref{7} uses a multilinear version of the Russo-Dye Theorem that we prove in this paper. This is of independent interest and is likely to be useful in other situations.

The norm of a multilinear map $\Phi$ from $\mathcal B(\Hil)^m$ into $\mathcal B(\K)$ is defined as
\begin{equation}
\|\Phi\|= \underset{\|X^{1}\| = \cdots = \|X^{m}\|=1}{\sup} \,\,\|\Phi (X^1, \ldots, X^m)\|. \label{9}
\end{equation}
We say $\Phi$ is \emph{positive} if $\Phi(X^1,\ldots,X^m)$ is a positive semidefinite operator whenever $X^1,\ldots,X^m$ are positive semidefinite. We prove:

\begin{theorem}\label{theorem2}
Let $\Phi$ be a positive multilinear map. Then
\begin{equation}
\|\Phi\|=\|\Phi(I,I,\ldots,I)\|. \label{10}
\end{equation}
\end{theorem}

The extensions of higher order analogues of Jacobi's formula \eqref{Jacobi} to derivatives of the antisymmetric tensor powers $\wedge^k (A)$ in \cite{tjain} were obtained by Jain by following an ``upwards from the bottom" approach, thinking of $\wedge^k (A)$ as an $\binom{n}{k} \times \binom{n}{k}$ matrix whose entries are $k\times k$ minors of $A.$ Grover \cite{grover} followed a similar approach in obtaining an expression for the derivatives of symmetric tensor powers $\vee^k(A), \ k\geq 1.$ Here we look at these problems following a ``downwards from the top" approach, thinking of $\wedge^k (A)$ and $\vee^k (A)$ as the restrictions of the tensor power $\otimes^k(A)$ to invariant subspaces. This has several advantages: the proofs become easier and more transparent, the formulas are seen to be valid for infinite dimensional operators as well, the path to studying the same problem for other symmetry classes of tensors becomes clearer.

\vskip0.1in
\vskip0.1in

For the most part, we concentrate on finite dimensional Hilbert spaces. Extensions to infinite dimensional spaces are briefly indicated. In Section 3 we provide the derivatives of the maps that take an operator to its $k$th tensor power, $k$th antisymmetric tensor power and $k$th symmetric tensor power and compute their norms. Closely related to the determinant is the permanent function. In Section 4 we give expressions for derivatives of all orders for the permanent. These supplement the results in \cite{rbhjain}. Results of this section have been reported in a survey article \cite{grover}. We provide the details here.

\vskip0.1in
\section{A Russo-Dye Theorem for Multilinear Maps}
We first provide a proof of Theorem \ref{theorem2}, when $\Hil$ is finite dimensional. We imitate the proof for positive linear maps given in \cite[p.41]{rBh}.
\vskip0.1in
Let $U^1, U^2, \ldots, U^m$ be unitary matrices and let
\begin{equation}
 U^i = \sum_{j=1}^{r_{i}} \,\,\lambda_j^i \,\, P_j^i, \quad 1 \leq i \leq m, \label{eq17}
\end{equation}
be their spectral resolutions. Here $\lambda_j^i$ are the eigenvalues of $U^i,$ and $P_j^i$ the corresponding eigenprojections. In particular, $|\lambda_j^i| = 1,$ $P_j^i$ are positive semidefinite, and
\begin{equation}
\sum_{j=1}^{r_{i}} P_j^i = I, \quad 1 \leq i \leq m.   \label{eq18}
\end{equation}
Since $\Phi$ is multilinear we have
\begin{eqnarray*}
\lefteqn{\Phi (U^1, \ldots, U^m)} \\
&=& \sum_{j_{1}=1}^{r_{1}} \sum_{j_{2}=1}^{r_{2}} \cdots \sum_{j_{m}=1}^{r_{m}} \lambda_{j_{1}}^1  \lambda_{j_{2}}^2 \ldots  \lambda_{j_{m}}^m \,\,\Phi \left (P_{j_{1}}^1, P_{j_{2}}^2, \ldots, P_{j_{m}}^m \right ),  
\end{eqnarray*}
and
$$\Phi (I, \ldots, I) =  \sum_{j_{1}=1}^{r_{1}} \sum_{j_{2}=1}^{r_{2}} \cdots \sum_{j_{m}=1}^{r_{m}}  \,\,\Phi \left (P_{j_{1}}^1, P_{j_{2}}^2, \ldots, P_{j_{m}}^m \right ).$$
Since $\Phi$ is positive, the operators $\Phi \left (P_{j_{1}}^1, P_{j_{2}}^2, \ldots, P_{j_{m}}^m \right )$ are positive semidefinite.
\vskip0.1in
Using these two relations, we see that
\begin{eqnarray*}
\lefteqn{\left [ \begin{array}{ll}\Phi (I,\ldots,I) & \Phi (U^1, \ldots, U^m) \\ \Phi (U^1, \ldots, U^m)^{\ast} & \Phi (I, \ldots, I) \end{array} \right ]} \\
&=& \sum_{j_{1}=1}^{r_{1}} \cdots \sum_{j_{m}=1}^{r_{m}} \left [ \begin{array}{cc} 1 & \lambda_{j_{1}}^1 \cdots  \lambda_{j_{m}}^m  \\ \overline{\lambda_{j_{1}}^1 \cdots  \lambda_{j_{m}}^m  } & 1 \end{array} \right ] \otimes \Phi  (P_{j_{1}}^1, \ldots, P_{j_{m}}^m).
\end{eqnarray*}\vskip0.1in
This is a sum of tensor products of positive semidefinite matrices, and is therefore, positive semidefinite. It follows from Proposition 1.3.2 in \cite{rBh} that
\begin{equation}
 \| \Phi (U^1, \ldots, U^m)\| \leq \| \Phi (I, \ldots, I)\|. \label{eq19}
\end{equation}
Now let $X^i,$ $1 \leq i \leq m,$ be matrices with $\|X^i\|=1.$ Then there exist unitary matrices $U^i$ and $V^i$ such that $X^i = \frac{1}{2} (U^i+V^i)$ (see \cite[p.42]{rBh}).
\vskip0.1in
By the multilinearity of $\Phi$
$$\Phi (X^1, \ldots, X^m) = \frac{1}{2^m} \sum \Phi (W^1, \ldots, W^m), $$
where the summation is over $2^m$ terms obtained by choosing each of the $W^i$ to be either $U^i$ or $V^i,$ $1 \leq i \leq m.$ It follows from \eqref{eq19} that
$$\| \Phi (X^1, \ldots, X^m)\| \leq \| \Phi (I, \ldots, I) \|. $$
Hence $\| \Phi\| = \| \Phi (I, \ldots, I)\|.$ This establishes Theorem \ref{theorem2} when $\Hil$ is finite dimensional.
\vskip0.2in
This theorem is also valid when $\Hil$ and $\K$ are infinite dimensional Hilbert spaces. It is likely to be useful, and we provide a proof for the infinite dimensional case. 

Our proof invokes the well-known fact that if $A$ and $B$ are positive operators on a Hilbert space, then $\left [ \begin{array}{cc} A & X \\ X^{\ast} & B \end{array} \right ]$ is positive if and only if there exists a contraction $K$ such that  $X = A^{1/2} \,\,K \,\,B^{1/2}.$ (See Theorem I.1 in \cite{ta}. This is Proposition 1.3.2 in \cite{rBh} but the proofs given there are only for finite dimensional spaces.) To prove \eqref{10} we have to show that if $X^1, \ldots, X^m$ are operators with $\|X^i\| \leq 1,$ then
\begin{equation}
 \|\Phi (X^1, \ldots, X^m)\| \leq \| \Phi (I, \ldots, I)\|.   \label{eq20}
\end{equation}
Consider first the case when
\begin{equation}
X^i = \sum_{j=1}^{r_{i}} \lambda_j^{i} \,\,P_j^i, \quad 1 \leq i \leq m,    \label{eq21}
\end{equation}
where $P_j^i$ are mutually orthogonal projection operators  with $\sum\limits_{j=1}^{r_{i}} \,\,P_j^i = I,$ and $|\lambda_j^i| = 1.$ It can be seen that $\|X^i\| \leq 1.$ (See \cite[p.11]{paulsen}.) Arguing as in the finite dimensional case we see that the inequality \eqref{eq20} holds in this case. Now if $U^i,$ $1 \leq i \leq m,$ are unitary operators, then each $U^i$ is a limit of a sequence of operators of the form \eqref{eq21}. This shows that the inequality \eqref{eq20} holds when $X^i$ are unitary. From here one can see that the inequality continues to hold if each $X^i$ is a convex combination of unitary operators. Finally, since the closed unit ball in $\mathcal{B}(\mathcal{H})$ is the closed convex hull of unitary operators (see \cite[p.75]{halmos}) the inequality is valid when $X^i$ are any operators with $\|X^i\| \leq 1.$
\vskip0.1in

\section{Formulas for $\De^m \otimes^k (A), \ \De^m \wedge^k (A), \ \De^m \vee^k (A)$ and their norms}
Let $\mathcal{H}$ be a Hilbert space and let $\otimes^k \mathcal{H}$  be its $k$-fold tensor power $\mathcal{H} \otimes \mathcal{H} \otimes \cdots \otimes \mathcal{H}.$ Let $\wedge^k \mathcal{H}$ and $\vee^k \mathcal{H}$ be the subspaces of $\otimes^k \mathcal{H}$  consisting of antisymmetric tensors and symmetric tensors, respectively. (See \cite[Chap. I]{Rbh} for definitions, notations and basic facts.) If $\dim \mathcal{H} = n,$ then $\dim \wedge^k \mathcal{H} =\binom{n}{k}$ for $1 \leq k \leq n,$ and $\dim \vee^k \mathcal{H}=\binom{n+k-1}{k}$ for $k\geq 1$. For $k > n,$ the space $\wedge^k  \mathcal{H}$ is taken to be zero. For every $A$ in $\mathcal{B}(\mathcal{H}),$ we denote by $\otimes^k A $ its $k$-fold tensor power $A \otimes A \otimes \cdots \otimes A.$ This is an operator on $\otimes^k \mathcal{H}$ that leaves invariant the subspaces $\wedge^k \mathcal{H}$ and $\vee^k\mathcal{H}.$ The restriction of $\otimes^k A$ to these subspaces are denoted by $\wedge^k A,$ the $k$th antisymmetric tensor power of $A$ and $\vee^k A$, the $k$th symmetric tensor power of $A$, respectively. We wish to describe the $m$th derivatives of the maps $A \mapsto \otimes^k A,$ $A \mapsto \wedge^k A,$ and $A \mapsto \vee^k A$, from $\mathcal{B}(\mathcal{H})$ into $\mathcal{B}(\otimes^k \mathcal{H}),\ \mathcal{B}(\wedge^k \mathcal{H}),$ and $\mathcal{B}(\vee^k \mathcal{H}),$ respectively. 

\vskip0.1in
If $f : X \rightarrow Y$ is a map between normed spaces, then its $m$th derivative at a point $a$ (if it exists) is a map $\De^m f (a)$ from the $m$-fold product $X \times \cdots \times X$ into $Y.$  As a function of its $m$ variables $\De^m f(a) (x^1, \ldots, x^m)$ is symmetric and linear in each variable. One way of computing it is by using the relation
\begin{eqnarray}
 \lefteqn{\De^m f(a) (x^1, \ldots, x^m)        } \nonumber \\
&=& \left . \frac{\partial^m}{\partial t_1  \cdots \partial t_m} \right |_{t_{1} = \cdots = t_{m} = 0} \,\, f(a + t_1 x^1 + \cdots + t_m x^m).\label{derdefn}
\end{eqnarray}
We refer the reader to Chapter X of \cite{Rbh} for basic facts about differential calculus on matrix spaces.
Given $A^1, \ldots, A^k$ in $\mathcal{B}(\mathcal{H})$ we define their {\it symmetrised tensor product} as 
\begin{equation}
 A^1 \widetilde{\otimes} A^2 \widetilde{\otimes} \cdots \widetilde{\otimes} A^k = \frac{1}{k!} \sum_{\sigma \in S_{k}} A^{\sigma(1)} \otimes A^{\sigma(2)} \otimes \cdots \otimes  A^{\sigma(k)}, \label{eq1} 
\end{equation}
where $S_k$ is the set of all permutations on $k$ symbols. One can check that this operator on $\otimes^k \mathcal{H}$ leaves invariant the subspaces $\wedge^k \mathcal{H}$ and ${\vee^k \mathcal{H}.}$ The restriction of the symmetrised tensor product to $\wedge^k \mathcal{H}$ and ${\vee^k \mathcal{H}}$ will be denoted by
\begin{equation}
A^1 \wedge A^2 \wedge \cdots \wedge A^k \text{ and } A^1\vee A^2 \vee \cdots \vee A^k,   \label{eq2}
\end{equation}
respectively and called the {\it symmetrised antisymmetric tensor product} and the {\it symmetrised symmetric tensor product} of $A^1, A^2, \ldots, A^k.$ The operator $A^1 \wedge A^2 \wedge \cdots \wedge A^k$ acts on product vectors $u_1 \wedge u_2 \wedge \cdots \wedge u_k$ as 
\begin{eqnarray}
\lefteqn{ (A^1 \wedge A^2 \wedge \cdots \wedge A^k)(u_1 \wedge u_2 \wedge \cdots \wedge u_k) }  \nonumber\\
&=& \frac{1}{k!} \sum_{\sigma \in S_{k}} (A^{\sigma(1)} u_1) \wedge (A^{\sigma(2)} u_2) \wedge \cdots \wedge (A^{\sigma(k)} u_k). \label{eq3}
\end{eqnarray}
Similarly the operator $A^1\vee A^2\vee \cdots \vee A^k$ acts on $u_1\vee u_2\vee \cdots u_k$ as 
\begin{eqnarray}
\lefteqn{ (A^1 \vee A^2 \vee \cdots \vee A^k)(u_1 \vee u_2 \vee \cdots \vee u_k) }  \nonumber\\
&=& \frac{1}{k!} \sum_{\sigma \in S_{k}} (A^{\sigma(1)} u_1) \vee (A^{\sigma(2)} u_2) \vee \cdots \vee (A^{\sigma(k)} u_k). \label{eq3.2}
\end{eqnarray}

Let $f$ be the real function $f(t) = t^k.$ Then $f^{(m)} (t) = k (k-1) \ldots (k-m+1) t^{k-m} = \frac{k!}{(k-m)!} \,t^{k-m} $ for $1 \leq m \leq k,$ and $f^{(m)} (t)=0$ for $m>k.$ The following theorem is an operator analogue of this. 
With the above notations, we have:

\begin{theorem} \label{thm1}
Let $1 \leq m \leq k.$ The $m$th derivatives of the maps $\otimes^k, \wedge^k$ and ${\vee^k}$ are given by the formulas
\begin{equation}
\De^m \otimes^k (A) (X^1, \ldots, X^m)
= \frac{k!}{(k-m)!}\ \underset{k-m\,\,copies}{\underbrace{A \widetilde{\otimes} \cdots \widetilde{\otimes} A }} \widetilde{\otimes} X^1 \widetilde{\otimes} X^2 \widetilde{\otimes} \cdots \widetilde{\otimes} X^m, \label{eq4}
\end{equation}
\begin{equation}
\De^m \wedge^k (A) (X^1, \ldots, X^m)= \frac{k!}{(k-m)!} \ \underset{k-m\,\,copies}{\underbrace{A \wedge \cdots \wedge A }}\wedge X^1 \wedge X^2 \wedge \cdots \wedge X^m,  \label{eq5}
\end{equation}
and
\begin{equation}
\De^m \vee^k (A) (X^1, \ldots, X^m) = \frac{k!}{(k-m)!}\  \underset{k-m\,\,copies}{\underbrace{A \vee \cdots \vee A }} \vee
X^1 \vee X^2 \vee \cdots \vee X^m.  \label{eq7}
\end{equation}
If $m > k,$ then all the derivatives are zero.
\end{theorem}
\vskip0.1in

\begin{proof}
By the formula \eqref{derdefn},
\begin{eqnarray}
 \lefteqn{\De^m \otimes^k (A) (X^1, \ldots, X^m)        } \nonumber \\
&=& \left . \frac{\partial^m}{\partial t_1  \cdots \partial t_m} \right |_{t_{1} = \cdots = t_{m} = 0} \,\, \otimes^k (A + t_1 X^1 + \cdots + t_m X^m). \label{eq11}
\end{eqnarray}

To evaluate this we expand the $k$-fold  tensor product on the right hand side. The resulting expansion is a polynomial in the variables $t_1, \ldots, t_m.$ The derivative in \eqref{eq11} is evidently the coefficient of the term $t_1 t_2 \cdots t_m$ in this polynomial. One can check that this is given by the expression \eqref{eq4}.
\vskip0.1in
Next we prove \eqref{eq5} using \eqref{eq4}. The proof for \eqref{eq7} is similar. 
The chain rule of differentiation for a composite function $g \circ f$ says that
$$\De (g \circ f) (a) (x) = \De g (f(a)) (\De f (a) (x)). $$
If $L$ is a linear map, then its derivative is equal to $L,$ and in this case
$$\De (L \circ f)(a) (x) = L (\De f (a) (x)). $$
Repeating this argument one sees that if $f$ is $m$ times differentiable, then
\begin{equation}
\De^m ( L \circ f)(a) (x^1, \ldots, x^m) = L (\De^m f (a) (x^1, \ldots, x^m)). \label{14.2}
\end{equation}

Now let $Q_k: \wedge^k \Hil \rightarrow \otimes^k \Hil$ be the inclusion map. Then $Q_k^*:\otimes^k \Hil \rightarrow \wedge^k \Hil$ is the projection given by
$$Q_k^*(x_1\otimes \cdots\otimes x_k)=\frac{1}{k!} \sum_{\sigma\in S_k} \varepsilon_{\sigma} x_{\sigma(1)}\otimes \cdots\otimes x_{\sigma(k)},$$
where $\varepsilon_{\sigma}=\pm 1,$ depending on whether $\sigma$ is an even or an odd permutation.
Define $\tilde{Q}_k: \lo(\otimes^k \Hil)\rightarrow \lo(\wedge^k \Hil)$ by
\begin{equation}
\tilde Q_k(T)=Q_k^* T Q_k.\label{defnofQ_k}
\end{equation}
Also $\wedge^k: \lo(\Hil)\rightarrow \lo(\wedge^k \Hil)$ factors through $\otimes^k : \lo(\Hil)\rightarrow \lo(\otimes^k \Hil)$ via $\tilde Q_k$ as $$\wedge^k A=\tilde Q_k(\otimes^k A) \quad \forall\ A\in \lo(\Hil).$$
Since $\tilde Q_k $ is linear, by \eqref{14.2} we have
\begin{equation}
\De^m \wedge^k (A)=\tilde Q_k \circ \De^m \otimes^k (A).\label{17.1}
\end{equation}
 Using this we obtain the expression \eqref{eq5}. 
\vskip0.1in

From \eqref{eq4} we see that
\begin{equation*}
 \De^k \otimes^k (A) (X^1, \ldots, X^k)= k!\,\, X^1 \widetilde{\otimes} X^2 \widetilde{\otimes} \cdots \widetilde{\otimes} X^k.
\end{equation*}
This expression does not involve $A.$ Hence $\De^m \otimes^k (A) = 0$ if $m>k.$ Similarly $\De^m \wedge^k (A) = 0$ and $\De^m \vee^k (A) = 0$ if $m>k.$
\end{proof}
\vskip0.1in
\begin{remark}
\rm By putting $k=n$ in \eqref{eq5}, we get $$ \De^m \det (A) (X^1, \ldots, X^m)=\frac{n!}{(n-m)!} \ \underset{n-m\,\,copies}{\underbrace{A \wedge \cdots \wedge A }}\wedge X^1 \wedge X^2 \wedge \cdots \wedge X^m.$$ This can also be written as $$\De^m \det (A) (X^1, \ldots, X^m)=\frac{n!}{(n-m)!} \Delta(A, \ldots, A, X^1,\ldots, X^m),$$ where $\Delta(A, \ldots, A, X^1,\ldots, X^m)$ denotes the mixed discriminant of the matrices \\$A, \ldots, A, X^1, \ldots, X^m.$ This is the same as Theorem 1 in \cite{rbhjain}.
\end{remark}

From these formulas we obtain the values of the norms of these derivatives. We separate the cases of $\|\De^m \otimes^k (A)\|$ and $\|\De^m \vee^k (A)\|$. The evaluation of these norms is independent of Theorem \ref{theorem2}, whereas we make essential use of this theorem in calculating $\|\De^m \wedge^k (A)\|.$

\begin{theorem} \label{thm2}
 For $1 \leq m \leq k,$ we have
\begin{equation}
\|\De^m \otimes^k (A) \| = \frac{k!}{(k-m)!} \,\, \|A\|^{k-m} \label{eq7.2}
\end{equation}
and
\begin{equation}
\|\De^m \vee^k (A) \| = \frac{k!}{(k-m)!} \,\, \|A\|^{k-m}. \label{eq7.3}
\end{equation}
\end{theorem}

\begin{proof}
To compute the norm $\|\De^m \otimes^k (A)\|,$ we first see that by definition of the symmetrised tensor product \eqref{eq1} and by the triangle inequality we get
\begin{equation*}
\| \underset{k-m\,\,copies}{\underbrace{A \widetilde{\otimes} \cdots \widetilde{\otimes} A }} \widetilde{\otimes} X^1 \widetilde{\otimes} X^2 \widetilde{\otimes} \cdots \widetilde{\otimes} X^m\|\leq \frac{1}{k!} \sum_{\sigma \in S_k}\| Y^{\sigma(1)} \otimes Y^{\sigma(2)}\otimes \cdots\otimes Y^{\sigma(k)}\|,
\end{equation*}
where $k-m$ of the $Y$'s are equal to $A$ and the rest are $X^1, X^2, \ldots, X^m.$
Each of the terms in the summation is equal to $\|A\|^{k-m} \|X^1\| \|X^2\| \cdots \|X^m\|.$
By the definition of the norm of a multilinear map \eqref{9} we obtain
$$\|\De^m \otimes^k (A) \|\leq \frac{k!}{(k-m)!} \|A\|^{k-m}.$$
Also note that 
$$\left\|\De^m \otimes^k (A) \left(\frac{A}{\|A\|}, \frac{A}{\|A\|}, \ldots, \frac{A}{\|A\|}\right)\right\|=\frac{k!}{(k-m)!} \|A\|^{k-m}.$$
This shows that
$$\|\De^m \otimes^k (A)\|\geq \frac{k!}{(k-m)!} \|A\|^{k-m}.$$
Hence we obtain \eqref{eq7.2}.
This argument works equally well in infinite dimensions.

Let $R_k: \vee^k \Hil \rightarrow \otimes^k \Hil$ be the inclusion map. 
Define $\tilde{R}_k: \lo(\otimes^k \Hil)\rightarrow \lo(\vee^k \Hil)$ by
\begin{equation}
\tilde R_k(T)=R_k^* T R_k.\label{defnofR_k}
\end{equation}
Then $\|\tilde R_k\|\leq 1$. Arguments similar to those in the proof of Theorem \ref{thm1} lead to an expression similar to \eqref{17.1}:
$$\De^m \vee^k (A)=\tilde R_k \circ \De^m \otimes^k (A).$$ It follows that
\begin{equation}
\|\De^m \vee^k (A)\|\leq \|\tilde R_k\| \|\De^m \otimes^k (A)\|\leq \frac{k!}{(k-m)!}\|A\|^{k-m}. \label{17}
\end{equation}

Let us now consider the case when $\mathcal{H}$ is an $n$ dimensional space. The polar decomposition theorem tells us that
$$A = U \,\,|A|,$$
where $U$ is unitary and $|A| = (A^{\ast} A)^{1/2}$ is positive semidefinite. Since $UU^{\ast} = I,$ we have
\begin{eqnarray*}
\lefteqn{\vee^k (A + t_1 X^1 + t_2 X^2 + \cdots + t_m X^m)       } \\
&=&  \vee^k \left (U (|A| + t_1 U^{\ast} X^1 + t_2 U^{\ast} X^2 + \cdots + t_m U^{\ast} X^m) \right )   \\
&=&(\vee^k  U)\, \vee^k  \left (|A| + t_1U^{\ast}X^1 + t_2U^{\ast} X^2 + \cdots + t_m U^{\ast} X^m \right ). 
\end{eqnarray*}
So from  \eqref{derdefn} we obtain
\begin{eqnarray}
\lefteqn{\De^m \vee^k (A) (X^1, \ldots, X^m)   } \nonumber \\
&=& (\vee^k U)\, \De^m \vee^k (|A|) (U^{\ast} X^1, \ldots, U^{\ast} X^m). \label{eq12}
\end{eqnarray}
Now $\vee^k U$ is unitary and the norm is unitarily invariant. So we have
\begin{eqnarray*}
\lefteqn{\| \De^m \vee^k (A) (X^1, \ldots, X^m) \|   } \\
&=& \| \De^m \vee^k (|A|) (U^{\ast} X^1, \ldots, U^{\ast} X^m) \|. 
\end{eqnarray*}
The condition $\| X^j\| = 1$ is equivalent to $\| U^{\ast} X^j\| = 1$ for $1 \leq j \leq m.$ So we have proved that
\begin{equation}
\| \De^m \vee^k (A) \| = \| \De^m \vee^k (|A|) \|.   \label{eq13}
\end{equation}

 Now assume $A$ is positive semidefinite and let $u$ be an eigenvector corresponding to its maximal eigenvalue $\|A\|.$ Consider the vector $w=u \vee u \vee \cdots \vee u$ in $\vee^k \mathcal{H}.$ If $T = Y^1 \vee Y^2\vee \cdots \vee Y^k$ is an operator in which $k-m$ of the $Y$'s are equal to $A$ and the rest of them are equal to  $I,$ then $Tw = \|A\|^{k-m} w.$ It then follows from \eqref{eq7} that
$$(\De^m \vee^k (A)(I,\ldots, I))w = \frac{k!}{(k-m)!} \,\,\|A\|^{k-m} w.$$ 
This shows that
$$\|\De^m \vee^k (A) \| \ge \frac{k!}{(k-m)!}  \|A\|^{k-m}.$$
We have already noted the reverse inequality in \eqref{17}. So we have \eqref{eq7.3} in the case when $A$ is positive semidefinite. The relation \eqref{eq13} then shows that \eqref{eq7.3} is valid for all $A$. 

We now indicate the modifications needed in this proof to handle the infinite dimensional case. In this case $A$ has a {\it maximal polar representation} $A = U|A|$ in which $U$ is either an isometry $(U^{\ast}U = I)$ or a coisometry $(UU^{\ast} = I)$ (\cite[p.75]{halmos}). When $\mathcal{H}$ is finite dimensional these two conditions are equivalent and $U$ is unitary. Our argument using the polar decomposition for proving \eqref{eq13} can be modified. A very similar idea is used in \cite{sunder} and we refer the reader to that paper for details.
\vskip0.1in
To prove \eqref{eq7.3} in the infinite dimensional case we may, therefore, again assume that $A$ is a positive operator. If $A$ has pure point spectrum, then the arguments given for the finite dimensional case serve equally well here. In particular, \eqref{eq7.3} is valid for compact operators. Every positive operator is a limit of a sequence of positive operators with pure point spectrum. Using this fact one can see that \eqref{eq7.3} is valid for all operators. 
\end{proof}

Note that for the above proof no use of Theorem \ref{theorem2} has been made. The formula \eqref{7} for $\|\De^m \wedge^k (A)\|$ is more interesting, and to prove it we do need to invoke Theorem \ref{theorem2}.
\vskip0.1in
To compute $\| \De^m \wedge^k (A) \|$ we first note that the symmetrised antisymmetric tensor product of positive semidefinite operators is positive semidefinite. It follows from \eqref{eq5} that if $A$ is positive semidefinite, then the map $\De^m \wedge^k (A)$ from $\mathcal{B}(\mathcal{H})^m$ into $\mathcal{B} (\wedge^k \mathcal{H})$ is a positive multilinear map. So, we have from Theorem \ref{theorem2}
\begin{equation}
 \| \De^m \wedge^k (A) \| = \| \De^m \wedge^k (A) (I, \ldots, I) \|. \label{eq16}
\end{equation}
Arguments similar to the ones used in the proof of Theorem \ref{thm2} show that
\begin{equation}
\| \De^m \wedge^k (A) \| = \| \De^m \wedge^k (|A|) \|.   \label{eq14}
\end{equation}
 So we assume $A$ to be positive semidefinite. By \eqref{eq16}, we have
\begin{eqnarray*}
\| \De^m \wedge^k(A) \| &=& \|\De^m \wedge^k (A) (I, \ldots, I) \| \\
&=&\left \| \left .\frac{\partial^m}{\partial t_1 \cdots \partial t_m} \right |_{t_{1} = \cdots = t_{m} = 0} \,\,\wedge^k (A + t_1 I + \cdots + t_m I)\right \|.
\end{eqnarray*}
By the spectral theorem there exists a unitary $W$ such that $A = WDW^{\ast},$ where $D$ is the diagonal matrix whose diagonal entries are $\alpha_1 \ge \cdots \ge \alpha_n(\ge 0),$ the eigenvalues of $A.$ The matrix $\wedge^k W$ is again unitary, and our norm is unitarily invariant. So in the right hand side of the equation above we can replace $A$ by $D.$ Now$$\wedge^k (D+ t_1 I + \cdots + t_m I)$$is a diagonal matrix of order $\left ( \begin{array}{c} n \\ k \end{array} \right ).$ Its norm is equal to its top diagonal entry, which is $$ \left .\frac{\partial^m}{\partial t_1 \cdots \partial t_m} \right |_{t_{1} = \cdots = t_{m} = 0} \,\, \prod_{j=1}^k (\alpha_j + t_1 + \cdots + t_m).$$
A calculation shows that this is equal to
$$m!\,\, p_{k-m} \,\,(\alpha_1, \ldots, \alpha_k).$$
This establishes \eqref{7} in the case when $A$ is positive semidefinite. The general case follows from \eqref{eq14}.
\vskip0.1in
Theorem \ref{theorem1} can be modified for infinite dimensional operators. The statement of this theorem involves the sequence $s_1 (A) \ge s_2(A) \ge \cdots .$ If we stretch the definitions and interpret a point of the essential spectrum of $|A|$ as an eigenvalue of infinite multiplicity, then Theorem \ref{theorem1} is valid for infinite dimensional operators too. The proof is similar to the proof for symmetric tensor powers.


\section{Formulas for $\De^m \per (A)$}\label{permanent}
 The \emph{permanent} of $A=(a_{ij})\in \mat$, written as $\per{A}$, is defined by
\begin{equation}
\per{A} = \sum_{\sigma \in S_n} a_{1 \sigma (1)}a_{2 \sigma (2)}\cdots a_{n \sigma (n)}.\label{defnofper}
\end{equation}
Since the definitions of $\per$ and $\det$ are similar, it is natural to expect a formula for $\De \per (A)$ similar to the Jacobi formula \eqref{Jacobi}. Applying the special case $m=1$ of \eqref{derdefn} to the $\per$ function, we see that $\De \per (A)(X)$ is the coefficient of $t$ in the polynomial ${\per (A+tX)}$.
 For $1\leq j\leq n$, let $A(j;X)$ be the matrix obtained from $A$ by replacing the $j^{th}$ column of $A$ by the $j^{th}$ column of $X$ and keeping the rest of the columns unchanged. Since $\per$ is a linear function in each of the columns, we get 
\begin{equation}
\De \per (A)(X)= \sum_{j=1}^{n}\per{A(j;X)}. \label{5}
\end{equation}
To give a formula analogous to the Jacobi formula, we define the \emph{permanental adjoint} of $A$ as the $n \times n$ matrix whose $(i,j)$-entry is $\per{A(i|j)}$, where $A(i|j)$ denotes the $ (n-1)\times(n-1)$ submatrix obtained from $A$ by deleting its $i$th row and $j$th column (see \cite[p.237]{merris}). Note that the adjugate of $A,$ $\adj A,$ is defined as the transpose of the matrix whose $(i,j)$-entry is $(-1)^{i+j} \det A(i|j)$, whereas in the definition of $\padj$, the transpose is not taken. This is just a matter of convention. The expression \eqref{5} can be rewritten as follows.

\begin{theorem}
For each $X \in \mat$,
\begin{equation}
\De\per{(A)}(X) = \tr(\padj(A)^t X).  \label{4} 
\end{equation}
\end{theorem}
\vskip0.1in
Our next aim is to obtain higher order derivatives of the permanent function. The expressions obtained are analogous to the ones for the $\det$ function given in \cite{rbhjain}. Applying \eqref{derdefn} to the $\per$ function, we see that $\mbox \De^m \per A(X^1,\ldots,X^m)$ is the coefficient of $t_1\cdots t_m$ in the expansion of $ \per(A+t_1 X^1+\cdots+t_m X^m)$. To write an explicit expression for this, we require some notations. 
\vskip0.1in
Let
$Q_{m,n}=\{(i_1,\ldots,i_m) |{ \ i_1,\ldots,i_m\in \N,}\ 1\leq i_1<\cdots<i_m\leq n\}.$ For $m>n,\, Q_{m,n}=\varnothing$ by convention. Let $G_{m,n}=\{(i_1,\ldots,i_m)|$ $ i_1,\ldots,i_m \in \N, 1\leq i_1\leq\cdots\leq i_m\leq n\}.$ 
Note that for $m\leq n,\ Q_{m,n}$ is a subset of $G_{m,n}$. For ${\J=(j_1,\ldots,j_m) \in Q_{m,n}}$, we denote by $A(\J;X^1,\ldots,X^m)$, the matrix obtained from $A$ by replacing the $j_p^{th}$ column of $A$ by the $j_p^{th}$ column of $X^p$ for $1\leq p \leq m$, and keeping the rest of the columns unchanged. Expanding $ \per(A+t_1 X^1+\cdots+t_m X^m)$ by using the fact that $\per$ is a linear function in each of the columns, we obtain an expression for $\De^m \per A(X^1,\ldots,X^m)$ as follows. This is a generalisation of \eqref{5}. 

\begin{theorem}\label{4.2}
For $1\leq m\leq n,$
\begin{equation}
\De^m \per{(A)(X^1,\ldots,X^m)}=\sum_{\sigma \in S_m} \sum_{\J \in Q_{m,n}} \per{A(\J;X^{\sigma(1)},X^{\sigma(2)},\ldots, X^{\sigma(m)})}.\label{8'}
\end{equation}
In particular, 
$$\De^m \per{(A) (X,\ldots,X)}=m! \sum_{\J \in Q_{m,n}} \per{A(\J;X,\ldots,X)}.$$ 
\end{theorem}
\vskip0.1in

The \emph{Laplace expansion theorem} for permanents \cite[p. 16]{minc} says that
for any $1\leq m\leq n,$ and for any $\I \in Q_{m,n}$,
\begin{equation}
\per{A}=\sum_{\J \in Q_{m,n}} \per{A[\I|\J]} \per{A(\I|\J)},\label{laplace}
\end{equation}
where $A[\I|\J]$ denotes the $m \times m$ submatrix obtained from $A$ by picking rows $\I$ and columns $\J$ and $A(\I|\J)$ denotes the ${(n-m)\times(n-m)}$ submatrix obtained from $A$ by deleting rows $\I$ and columns $\J.$
In particular, for any $i, 1\leq i\leq n,$
\begin{equation}
\per{A}=\sum_{j=1}^{n} a_{ij}\, \per{A(i|j)}.\label{2'}
\end{equation} 
Using this, equation \eqref{5} can be rewritten as
\begin{equation}
\De\per{(A)}(X)=\sum_{i=1}^{n} \sum_{j=1}^{n} x_{ij} \per{A(i|j)}. \label{6}
\end{equation}
We obtain a generalisation of this expression for higher order derivatives. Let $Y_{[j]}$ denote the $j$th column of the matrix $Y$. Let $\sigma$ be a permutation on $m$ symbols, then by $Y^{\sigma}_{[\J]}$, we mean the matrix in which $Y^{\sigma}_{[j_p]}=X^{\sigma(p)}_{[j_p]}$ for $1\leq p \leq m$ and $Y^{\sigma}_{[\ell]}=0$ if $\ell$ does not occur in $\J.$ By using the Laplace expansion \eqref{laplace} for each term in the summation of \eqref{8'}, we obtain the following expression for $\De^m \per{(A)(X^1,\ldots,X^m)}$.

\begin{theorem}\label{thm7}
For $1\leq m\leq n,$
\begin{equation}
\De^m \per{(A)(X^1,\ldots,X^m)}=\sum_{\sigma \in S_m} \sum_{\I,\J \in Q_{m,n}} \per{A(\I|\J)} \per{Y^{\sigma}_{[\J]}[\I|\J]}. \label{2.9}
\end{equation}
In particular, 
$$\De^m \per{(A)(X,\ldots,X)}=m! \sum_{\I,\J \in Q_{m,n}} \per{A(\I|\J)} \per{X[\I|\J]}.$$
\end{theorem} 
\vskip0.1cm
Note that 
\begin{equation}
\De^n \per{(A)(X,\ldots,X)}=n!\ \per{X},
\end{equation}
and 
\begin{equation}
\De^m \per{(A)(X^1,\ldots,X^m)}=0 \text{ for all } m > n.
\end{equation}

We now describe a generalisation of \eqref{4} for higher order derivatives of the $\per$ function. 
Given an orthonormal basis $\{e_1,e_2,\ldots,e_n\}$ of an $n$ dimensional Hibert space $\Hil,$ 
the set $\{m(\alpha)^{-1/2}\ e_{\alpha}: \alpha \in G_{m,n}\}$ is an orthonormal basis of $\vee^m \Hil$. (See \cite[p.17]{Rbh} for details.)
Let $P_m$ be the canonical projection of $\vee^m \Hil$ onto the subspace ${\{e_{\alpha} : \alpha \in Q_{m,n}\}}$. Then there is a permutation of the above orthonormal basis of $\vee^m \Hil$ in which 
$P_m=\left[\begin{array}{ccc} 
I &\ O \\
O &\ O
\end{array}\right]$
and the matrix $T_m,$ defined by $T_m=\left(\per A[\alpha|\beta]\right)_{\alpha,\beta\in Q_{m,n}},$ is the upper left corner of $\vee^m A,$ that is, 
$$P_m \left(\vee^m A\right) P_m=
\left[\begin{array}{cc} 
T_m & \ O\\
O  & \ O
\end{array}\right].$$
Let $U$ be the $\binom{n}{m} \times \binom{n}{m}$ unitary matrix given by
$U=\left[\begin{array}{cccc}
\ & \ & \ &\  1\\
\ & \ &\ 1 & \ \\
\ & \ \text{\rotatebox{90}{\mbox{$\ddots$}}}& \ & \ \\
1& \ & \ & \ 
\end{array}\right].$
 For $\alpha,\beta \in Q_{n-m,n},$ the $(\alpha,\beta)$-entry of $U^* T_m U$ is $\per A(\alpha|\beta)$. Let $\widetilde{U}$ be the $\binom{n+m-1}{m}\times \binom{n+m-1}{m}$ matrix given by 
$\widetilde{U} =
\left[\begin{array}{cc} 
U & \ O\\
O & \  I
\end{array}\right].$
Let $\widetilde {\vee}^m A$ denote the matrix $\widetilde{U}^* \left(\vee^m A\right)^t \widetilde{U}$. Then 
\begin{equation}
P_m \left(\widetilde {\vee} ^m A \right)P_m=\left[\begin{array}{cc}
U^* T_m^t U & \ O\\
O & \ O
\end{array}\right].
\end{equation}
In particular for $m=n-1,$
$$P_{n-1} \left(\widetilde {\vee}^{n-1} A\right)P_{n-1}=\left[\begin{array}{ccc}
(\padj A)^t & \ O\\
O  & \ O
\end{array}\right].$$
Identifying an $n\times n$ matrix $X$ with $\binom{2n-2}{n-1}\times \binom{2n-2}{n-1}$ matrix
$\left[\begin{array}{ccc} 
X & \ O\\
O & \ O\end{array}\right]$, equation \eqref{4} can be rewritten as 
\begin{equation}
\De\per{(A)(X)}=\tr{\left(P_{n-1}\left(\widetilde{\vee}^{n-1} A\right)P_{n-1}\right)X}.\label{eq10}
\end{equation}
Its generalisation for higher order derivatives is given as follows.

\begin{theorem}\label{Theorem 9}
For $1\leq m\leq n$,
\begin{eqnarray}
\lefteqn{\De^m \per{(A)(X^1,\ldots,X^m)}}\nonumber\\
&=&m!\ \tr\left[\left(P_{n-m}\left(\widetilde{\vee}^{n-m} A \right)P_{n-m}\right)\left(P_m\left(X^1\vee\cdots\vee X^m \right)P_m\right)\right].\nonumber \\
&& \label{2.17}
\end{eqnarray}
In particular,
$$\De^m \per{(A)(X,\ldots,X)}=m!\ \tr{\left[\left(P_{n-m}(\widetilde{\vee}^{n-m}A)P_{n-m}\right)\left(P_m \left(\vee^m X \right)P_m\right)\right]}.$$
\end{theorem}
\vskip0.1cm
To see a proof of Theorem \ref{Theorem 9}, we first describe the notion of \emph{mixed permanent} of $m \times m$ matrices $T^1,\ldots,T^m$. (This was first introduced by Bapat in \cite{bapat}.) It is denoted by $\Delta_p(T^1,\ldots,T^m)$, and is defined as 
$$\Delta_p(T^1,\ldots,T^m)=\frac{1}{m!}\sum_{\sigma \in S_m} \per \left[T^{\sigma(1)}_{[1]},\ldots,T^{\sigma(m)}_{[m]}\right].$$
When all $T^j=T$, then $\Delta_p(T,\ldots,T)=\per T.$
Observe that for $\I,\J \in Q_{m,n}$,
\begin{equation}
\sum_{\sigma\in S_m} \per Y^{\sigma}_{[\J]}[\I|\J]=m!\  \Delta_p(X^1[\I|\J],\ldots, X^m[\I|\J]).\label{observation2}
\end{equation}
Using this, Theorem \ref{thm7} can be rewritten as follows:
\begin{eqnarray}
\lefteqn{\De^m \per{(A)(X^1,\ldots,X^m)}}\label{50}\\ 
&=&m!\ \sum_{\I,\J \in Q_{m,n}} \per{A(\I|\J)}\ \Delta_p(X^1[\I|\J],\ldots, X^m[\I|\J]). \nonumber
\end{eqnarray}

Next we note that for $\I,\J \in G_{m,n}$, the $(\I,\J)$-entry of ${ X^1\vee\cdots\vee X^m}$ is 
\begin{equation}
(m(\I) m(\J))^{-1/2} \Delta_p(X^1[\I|\J],\ldots,X^m[\I|\J]).\label{observation1}
\end{equation}

In particular, if $\I,\J \in Q_{m,n}$, then the $(\I,\J)$-entry of ${ X^1\vee\cdots\vee X^m}$ is \\${\Delta_p(X^1[\I|\J],\ldots,X^m[\I|\J])}$. The $(\J,\I)$-entry of $P_{n-m}(\widetilde{\vee}^{n-m} A)P_{n-m}$ is ${\per A(\I|\J)}$. The expression \eqref{2.17} can now be easily seen as a reformulation of \eqref{50}. 

\section{Remarks}
\begin{enumerate}
\item[(1)] An alternative proof of Theorem \ref{4.2} can be given using \eqref{eq7}. We know that $\De^m \per (A)(X^1,\ldots,X^m)$ is the $(\alpha,\alpha)$-entry of ${\De^m\vee^n A(X^1,\ldots,X^m)}$ for $\alpha=(1,\ldots,n)$, which by \eqref{eq7} and \eqref{observation1} is $\frac{n!}{(n-m)!} \Delta_p(X^1,\ldots,X^m,A,\ldots,A).$ This is the same as $$\frac{1}{(n-m)!} \sum_{\tau \in S_n} \per \left[Y^{\tau(1)}_{[1]},\ldots,Y^{\tau(n)}_{[n]}\right],$$ where $k-m$ of the $Y$'s are equal to $A$ and the rest are $X^1, X^2, \ldots, X^m.$ 
For any given ${\J\in Q_{m,n}}$ and $\sigma \in S_m$, there are $(n-m)!$ terms in this summation which are equal to $\per{A(\J;X^{\sigma(1)},X^{\sigma(2)},\ldots, X^{\sigma(m)})}$. This gives Theorem \ref{4.2}.
\vskip0.1in
\item[(2)] An upper bound for norms of the derivatives of the permanent function can be obtained by using \eqref{eq7.3}. By using the fact that  $\De^m \per (A)(X^1,\ldots,X^m)$ is one of the entries of the matrix ${\De^m\vee^n A(X^1,\ldots,X^m),}$ we obtain
\begin{equation}
\|\De^m \per (A)\|\leq \frac{n!}{(n-m)!}\|A\|^{n-m}.\label{52}
\end{equation}
While we have equality in \eqref{eq7.3}, we may have strict inequality here. For example, let $A=\left[\begin{array}{ccc} 1 & \ \\ \ & 0 \end{array}\right]$. Then $A$ is a positive semidefinite matrix. So $\De \per (A)$ is a positive linear functional. By the Russo-Dye Theorem, we have
\begin{equation}
\|\De \per (A)\|=|\De \per (A)(I)|,
\end{equation}
which is equal to 1, by \eqref{6}. But the right hand side of \eqref{52} is equal to 2.
\vskip0.1in
\item[(3)] In this paper we have limited ourselves to tensor powers, symmetric tensor powers and antisymmetric tensor powers. There are other symmetry classes of tensors, and the corresponding problems for these classes have been studied by Carvalho and Freitas in \cite{carvalhofreitas1} and \cite{carvalhofreitas2}. Norms of first derivatives of    the operators induced on the symmetry classes of tensors had been computed earlier by Bhatia and Da Silva \cite{rbhsilva}. The work in \cite{carvalhofreitas2} extends this to higher order derivatives.
\end{enumerate}

\end{document}